	\documentclass[12pt, a4paper]{article}
	
 	\usepackage{latexsym}
	\usepackage{amsmath}	
	\usepackage{epsfig}	
	\usepackage{times}	
	\usepackage{amssymb}
	
 	\newtheorem{thr}{Theorem}[section]
 	\newtheorem{prop}{Proposition}[section]

	\newtheorem{prob}{Problem}[section]
\begin{document}
 	\centerline{\Large{\bf On the number of permutation polynomials}}
 	\centerline{}
	\centerline{\Large{\bf over a finite field}}
	\centerline{}
 	\centerline{\textsuperscript{a}Kwang-Yon Kim, \textsuperscript{b}Ryul Kim}
 	\centerline{}
 	{\small \centerline{Faculty of Mathematics,  \textbf{Kim Il Sung} University, 
	Pyongyang, D.P.R Korea}}
	{\small \centerline{{e-mail address : \textsuperscript{a}kimkwangyon@yahoo.com, \textsuperscript{b}ryul\_ kim@yahoo.com}}}
	\centerline{}
	\centerline{}

 %
	\begin{abstract}
	We find a formula for the number of permutation polynomials of degree $q-2$ over a finite field $\textbf{F}_q$, 
	which has $q$ elements, in terms of the permanent of a matrix. 
	We write down an expression for the number of permutation polynomials of degree $q-2$ over a finite field $\textbf{F}_q$, 
	using the permanent of a matrix whose entries are $p$th roots of unity and using this obtain a nontrivial bound for the number. 
	Finally, we provide a formula for the number of permutation polynomials of degree $d$ less than $q-2$. 
	\end{abstract}
	{\bf Keywords} \small {Finite field, Permutation polynomial, Permanent}\\

%
%
%
	\section{Introduction}

	Let $\textbf{F}_q$ be the finite field of $q=p^r$ elements, where $p$ is a prime number and $r$ is a positive integer.
	A polynomial $f$ over $\textbf{F}_q$ is said to be a $\textit{permutation polynomial}$ if the induced map $\alpha \mapsto f(\alpha)$ 
	from $\textbf{F}_q$ to itself is a permutation. Let $N_q(d)$  denote the number of all permutation polynomials over $\textbf{F}_q$ 
	such that $\textnormal{deg}(f)=d$ and $f(0)=0$.
	There has been tremendous amount of interest in permutation polynomials over $\textbf{F}_q$ due to their applications 
	in cryptography, coding theory and combinatorics. 

	The main topics on the permutation polynomials over finite fields 
	are concerned with the problems of finding new classes of permutation polynomials and determining the number of 
	permutation polynomials with restriction for their degree or coefficients. 
	Recently there has been significant progress in finding new classes of permutation polynomials; 
	for instance, see \cite{akb,fer,hou,lai,lih,mao,sha,tuz,wul,zha}.

	In [9], Lidl and Mullen proposed nine open problems and conjectures involving permutation polynomials of finite fields. 
	This is one of the open problems yet.	
	
%
%
	\begin{prob} 
	(Lidl-Mullen). Let $N_q(d)$ denote the number of permutation polynomials of $\textnormal{\textbf{F}}_q$ which have degree $d$. 
	We have the trivial boundary conditions: $N_q(1)=q(q-1), N_q(d)=0$ if $d$ is a divisor of $q-1$ larger than 1, 
	and $\sum N_q(d)=q!$ where the sum is over all $1 \leq d<q-1$ such that $d$ is either 1 or it is not a divisor 
	of $q-1$. Find $N_q(d)$.
	\end{prob}

	In \cite{das}, Pinaki Das related the number of permutation polynomials in $\textbf{F}_q[x]$ of degree $d \leq q-2$ 
	to the solutions $(x_1, x_2, \cdots, x_q)$ of a system of linear equations over $\textbf{F}_q$, with the added restriction 
	that $x_i \neq 0$ and $x_i \neq x_j$ for $i \neq j$.

	Let $\omega \in \textbf{F}_q^{\times}$ be a primitive element of $\textbf{F}_q$.

%
%
	\begin{thr} 
	\textnormal{\cite{das}}
	\begin{equation}
	N_q(q-2)=(q-1)!-\# \big( x_1+\omega x_2 +\omega^2x_3+\cdots+\omega^{q-2}x_{q-1}=0 \big),
	\end{equation}
	where \# denotes the number of solutions in $\textnormal{\textbf{F}}_q$ of the above equation with the restriction 
	that $x_i \neq 0$ and $x_i \neq x_j$ for $i \neq j$.
	\end{thr}

%
%
	\begin{thr} 
	\textnormal{\cite{das}} Let $G_q(d)$ be the number of solutions in $\textnormal{\textbf{F}}_q$ of the system of equations
	\begin{align*}
	& x_1+\omega^{q-d-1} x_2 +\omega^{2(q-d-1)}x_3+\cdots+\omega^{(q-2)(q-d-1)}x_{q-1}=0 \\
	& x_1+\omega^{q-d-2} x_2 +\omega^{2(q-d-2)}x_3+\cdots+\omega^{(q-2)(q-d-2)}x_{q-1}=0 \\
	& \cdots \quad \cdots \quad \cdots \quad \cdots \quad \cdots \\
	& x_1+\omega x_2 +\omega^2x_3+\cdots+\omega^{q-2}x_{q-1}=0,
	\end{align*}
	such that $x_i \neq 0$ and $x_i \neq x_j$ for $i \neq j$. Then we have
	\begin{equation*}
	N_q(d)=(q-1)!-N_q(q-2)-N_q(q-3)- \cdots -N_q(d+1)-G_q(d).
	\end{equation*}
	\end{thr}
	
	Pinaki Das found a formula for the number of permutation polynomials of degree $p-2$ over $\textbf{F}_p$, 
	using permanent of a matrix.
	
	For an $n \times n$ matrix $A= \big( a_{ij} \big)$, a $\textit{permanent}$ of $A$ is defined by
	\begin{equation*}
	per(A)=\sum_{\sigma \in S_n} \prod_{i=1}^n a_{i \sigma (i)}.
	\end{equation*}
	(See \cite{mau}).

%
%
	\begin{thr} 
	\textnormal{\cite{das}} Let $A=Vandermonde \big( X, \cdots, X^{p-1} \big)$ be the matrix defined by 
	$\big( X^{(i-1)j} \big) _{i, j=\overline{1, p-1}}$ and also let $per(A)=\sum c_i X^i$. Then
	\begin{equation*}
	\# \big( x_1+2x_2+3x_3+ \cdots +(p-1)x_{p-1}\equiv 0 \pmod{p} \big) = \sum_{i: p|i}c_i.
	\end{equation*}
	Hence
	\begin{equation*}
	N_p(p-2)=(p-1)!-\sum_{i: p|i}c_i.
	\end{equation*}
	\end{thr}
%
%
	\begin{thr} 
	\textnormal{\cite{das}} Let	
	\begin{equation*}
	A=Vandermonde \big( z_1z_2 \cdots z_n, z_1^2z_2^{2^2} \cdots z_n^{2^n}, \cdots, z_1^{p-1}z_2^{(p-1)^2} 
	\cdots z_n^{(p-1)^n} \big)
	\end{equation*}
	where $1 \leq n \leq p-2$. Also let $per(A)=\sum c_{i_1i_2 \cdots i_n} z_1^{i_1}z_2^{i_2} \cdots z_n^{i_n}$. 
	Then the number of solutions in $\textnormal{\textbf{F}}_p$ of the system of equations
	\begin{align*}
	& x_1+2^n x_2 +3^n x_3+\cdots+(p-1)^n x_{p-1}=0 \\
	& x_1+2^{n-1} x_2 +3^{n-1} x_3+\cdots+(p-1)^{n-1} x_{p-1}=0 \\
	& \cdots \quad \cdots \quad \cdots \quad \cdots \quad \cdots \\
	& x_1+2x_2 +3x_3+\cdots+(p-1)x_{p-1}=0,
	\end{align*}
	such that $x_i \neq 0$ and $x_i \neq x_j$ for $i \neq j$ is equal to
	\begin{equation*}
	\sum_{p|i_1, \cdots, p|i_n} c_{i_1 i_2 \cdots i_n},
	\end{equation*}
	where the sum is over all those coefficients $c_{i_1 i_2 \cdots i_n}$ for which p divides the exponent $i_k$ of each $z_k$.
	\end{thr}

	Also in \cite{das}, Pinaki Das derived an expression for the number of permutation polynomials of degree  $p-2$ 
	over $\textbf{F}_p$ in terms of the permanent of a Vandermonde matrix whose entries are the primitive $p$th roots of unity.

%
%
	\begin{thr} 
	\textnormal{\cite{das}} Let	$\zeta = e^{2\pi i/p}$ be a primitive $p$th root of unity and
 	\begin{equation*}
	V=Vandermonde \big( \zeta, \zeta^2, \cdots, \zeta^{p-1} \big).
	\end{equation*}
 Then
	\begin{equation*}
	N_p(p-2)=\Big( 1-\frac{1}{p} \Big) \big( (p-1)!-per(V) \big).
	\end{equation*}
	\end{thr}

%
%
	\begin{thr} 
	\textnormal{\cite{das}} Let	$V=Vandermonde \big( \zeta, \zeta^2, \cdots, \zeta^{p-1} \big)$, where $\zeta$ is a primitive pth root of unity. 
	Then
 	\begin{equation*}
	|per(V)| \leq \sqrt{\frac{1+(p-2)p^{p-1}}{p-1}}.
	\end{equation*}
 	\end{thr}
	So the following bounds are obtained.

%
%
	\begin{thr} 
	\textnormal{\cite{das}} 
 	\begin{equation*}
	\Bigg \arrowvert N_p(p-2)- \Big( 1-\frac{1}{p} \Big) (p-1)! \Bigg \arrowvert \leq \Big( 1-\frac{1}{p} \Big) \sqrt{\frac{1+(p-2)p^{p-1}}{p-1}}.
	\end{equation*}
 	\end{thr}

	Koyagin and Pappalardi proved the asymptotic formulas for the number of permutations for which the 
	associated permutation polynomial has some restriction for the degree or the coefficients(see \cite{ko1,ko2}). 
	We generalize the results of \cite{das} using the concepts of so-called formal polynomial and permanent.

%
%
 %
%
	\section{The number of permutation polynomials of degree $q-2$ over $\textbf{F}_q$}

	Let $\textbf{Z} \textbf{F}_q[X]$ be the set of all functions from $\textbf{F}_q$ to the set of integers $\textbf{Z}$ with three operations defined 
	by the following: for any $f, g \in \textbf{Z} \textbf{F}_q[X]$ and any $a \in \textbf{Z}$,
	\begin{equation*}
	f+g : \textbf{F}_q \to \textbf{Z}, ~ x \mapsto f(x)+g(x) 
	\end{equation*}
	\begin{equation*}
	f \cdot g : \textbf{F}_q \to \textbf{Z}, ~ x \mapsto \sum_{u+v=x}f(u) \cdot g(v) 
	\end{equation*}
	\begin{equation*}
	a \cdot f : \textbf{F}_q \to \textbf{Z}, ~ x \mapsto af(x).
	\end{equation*}
	Given an element $\alpha$ in $\textbf{F}_q$, $X^{\alpha}$ denotes the function from $\textbf{F}_q$ to $\textbf{Z}$ such that
	\begin{equation*}
	X^{\alpha}(x)=\left \{ 
	\begin{array}{ll}
		1, & x= \alpha \\
		0, & x \neq \alpha
	\end{array} \right..
	\end{equation*}

	Let $M=\big\{ X^{\alpha} \big \arrowvert \alpha \in \textbf{F}_q \big \}$. Then we can conclude the following.
%
%
	\begin{prop} 
	$\textnormal{\textbf{Z}\textbf{F}}_q[X]$ is a commutative free algebra over $\textnormal{\textbf{Z}}$ with a basis $M$ under the above three operations.
	\end{prop}

	Elements of this algebra are called $\it{formal \- polynomials}$.

%
%

	\begin{thr} 
	Let $A$ be the $(q-1) \times (q-1)$ matrix defined by
	\begin{equation*}
	\left(
	\begin{array}{cccc}
	X^{\omega} & X^{\omega^2} & \cdots & X^{\omega^{q-1}}\\
	X^{\omega^2} & X^{\omega^3} & \cdots & X^{\omega}\\
	\cdots & \cdots & \cdots & \cdots \\
	X^{\omega^{q-1}} & X^{\omega} & \cdots & X^{\omega^{q-2}}
	\end{array} \right)
	\end{equation*}
	and $per(A)=\sum_{i=0}^{q-2}c_i X^{\omega^i}+c_{-1}$. Then 
	\begin{equation*}
	\# \big( x_1+\omega x_2+\omega^2 x_3+ \cdots +\omega^{q-2} x_{q-1}=0 \big) = c_{-1}.
	\end{equation*}
	Hence
	\begin{equation}
	N_q(q-2)=(q-1)! -c_{-1}
	\end{equation}
	and also
	\begin{equation}
	c_{-1}+(q-1)c_0=(q-1)!.
	\end{equation}
	\end{thr}
	
	$\it{Proof}$. Since
	\begin{equation*}
	A=\left(
	\begin{array}{cccc}
	X^{\omega} & X^{\omega^2} & \cdots & X^{\omega^{q-1}}\\
	X^{\omega^2} & X^{\omega^3} & \cdots & X^{\omega}\\
	\cdots & \cdots & \cdots & \cdots \\
	X^{\omega^{q-1}} & X^{\omega} & \cdots & X^{\omega^{q-2}}
	\end{array} \right)
	=\left(
	\begin{array}{cccc}
	X^{\omega} & X^{\omega^2} & \cdots & X^{\omega^{q-1}}\\
	X^{\omega^{1+1}} & X^{\omega^{2+1}} & \cdots & X^{\omega^{q-1+1}}\\
	\cdots & \cdots & \cdots & \cdots \\
	X^{\omega^{1+q-2}} & X^{\omega^{2+q-2}} & \cdots & X^{\omega^{q-1+q-2}}
	\end{array} \right),
	\end{equation*}
	every summand of $per(A)$ is
	\begin{equation*}
	X^{\omega^{y_1}} X^{\omega^{1+y_2}} X^{\omega^{2+y_3}} \cdots  X^{\omega^{q-2+y_{q-1}}} = 
	X^{\omega^{y_1}+\omega \omega^{y_2}+\omega^2 \omega^{y_3}+\cdots +\omega^{q-2} \omega^{y_{q-1}}},
	\end{equation*}
	where $1 \leq y_i \leq q-1, y_i \neq y_j (i \neq j)$. So if we set $x_i=\omega^{y_i} (i=\overline{1, q-1})$, then
	\begin{equation*}
	per(A)=\sum_{\substack{x_i \in \textnormal{\textbf{F}}_q^{\times} \\ x_i \neq x_j ~ (i \neq j)}}X^{x_1+\omega x_2+
	\omega^2 x_3+\cdots +\omega^{q-2} x_{q-1}}.	
	\end{equation*}
	Hence $c_{-1}$ is equal to the number of solutions in $\textbf{F}_q^\times$ of $x_1+\omega x_2+\omega^2 x_3+\cdots +
	\omega^{q-2} x_{q-1}=0$ with the restriction that $x_i \neq x_j$ for $i \neq j$. And then we have $N_q(q-2)=(q-1)!-c_{-1}$ by (1).

	If $x_1+\omega x_2+\omega^2 x_3+\cdots +\omega^{q-2} x_{q-1}=1, x_i \neq 0$ and $x_i \neq x_j $ for $i \neq j$, 
	then for any integer $k$ and $x_i^{\prime}=\omega^k x_i, ~ x_1^{\prime}+\omega x_2^{\prime}+\omega^2 x_3^{\prime}+\cdots
	+\omega^{q-2}x_{q-1}^{\prime}=\omega^k$, with $x_i^{\prime} \neq 0$ and $x_i^{\prime} \neq x_j^{\prime}$ for $i \neq j$.

	Conversely, if $x_1+\omega x_2+\omega^2 x_3+\cdots +\omega^{q-2} x_{q-1}=\omega^k, x_i \neq 0$ and $x_i \neq x_j $ for $i \neq j$, 
	then setting $x_i^{\prime}=\omega^{-k} x_i$  gives that $x_1^{\prime}+\omega x_2^{\prime}+\omega^2 x_3^{\prime}+\cdots
	+\omega^{q-2}x_{q-1}^{\prime}=1$, with  $x_i^{\prime} \neq 0$ and $x_i^{\prime} \neq x_j^{\prime}$ for $i \neq j$. 
	Hence $c_0=c_k$ for all integers $k \in \{ 1, 2, \cdots, q-2 \}$. Since $\sum_{i=0}^{q-2}c_i+c_{-1}=(q-1)!$, 
	we have $c_{-1}+(q-1)c_0=(q-1)!$. $\Box$\\
	
	$\textbf{Example 2.1}$. The only permutation polynomials of degree 2 over $\textbf{F}_{2^2}$ are three monomials 
	$\omega X^2, \omega^2 X^2, \omega^3 X^2=X^2$, where $\omega$ is a primitive element of $\textbf{F}_{2^2}$.

	On the other hand, the permanent of the matrix
	\begin{equation*}
	A=\left(
	\begin{array}{ccc}
	X^{\omega} & X^{\omega^2} & X^{\omega^3}\\
	X^{\omega^2} & X^{\omega^3} & X^{\omega}\\
	X^{\omega^3} & X^{\omega} & X^{\omega^2}
	\end{array} \right)
	\end{equation*}
	is $per(A)=X^{\omega^3}+X^{\omega^2}+X^{\omega}+3$ and $c_{-1}=3$, thus $3!-c_{-1}=N_{2^2}(2)=3$.\\

	$\textbf{Example 2.2}$. When $\omega$ is a primitive element of $\textbf{F}_{2^3}$, the permanent of the matrix
	\begin{equation*}
	A=\left(
	\begin{array}{ccccccc}
	X^{\omega} & X^{\omega^2} & X^{\omega^3} & X^{\omega^4 }& X^{\omega^5} & X^{\omega^6} & X^{\omega^7} \\
	X^{\omega^2} & X^{\omega^3} & X^{\omega^4} & X^{\omega^5} & X^{\omega^6} & X^{\omega^7} & X^{\omega} \\
	X^{\omega^3} & X^{\omega^4} & X^{\omega^5} & X^{\omega^6} & X^{\omega^7} & X^{\omega} & X^{\omega^2} \\
	X^{\omega^4} & X^{\omega^5} & X^{\omega^6} & X^{\omega^7} & X^{\omega} & X^{\omega^2} & X^{\omega^3} \\
	X^{\omega^5} & X^{\omega^6} & X^{\omega^7} & X^{\omega} & X^{\omega^2} & X^{\omega^3} & X^{\omega^4} \\
	X^{\omega^6} & X^{\omega^7} & X^{\omega} & X^{\omega^2} & X^{\omega^3} & X^{\omega^4} & X^{\omega^5} \\
	X^{\omega^7} & X^{\omega} & X^{\omega^2} & X^{\omega^3} & X^{\omega^4} & X^{\omega^5} & X^{\omega^6} 
	\end{array} \right)
	\end{equation*}
	is
	\begin{align*}
	per(A) & = 624X^{\omega^7}+624X^{\omega^6}+624X^{\omega^5}+624X^{\omega^4}\\
			& + 624X^{\omega^3}+624X^{\omega^2}+624X^{\omega}+672,
	\end{align*}
	and $c_{-1}=672$. Thus $N_{2^3}(6)=7!-C_{-1}=5040-672=4368$. In fact, we can see by counting that there are $4368$ permutation 
	polynomials of degree $6$ over $\textbf{F}_{2^3}$.

	When
 	\begin{equation*}
	f(X)=\sum_{i=0}^{q-2}a_i X^{\omega^i}+a_{-1}
	\end{equation*}
	is a formal polynomial and $z$ is a complex number, we denote
 	\begin{equation*}
	\sum_{i=0}^{q-2}a_i z^{Tr(\omega^i)}+a_{-1}
	\end{equation*}
	as $f(z)$, where $Tr$ is the absolute trace in $\textbf{F}_q$.
%
%
	\begin{thr}
	Let $\zeta=e^{2\pi i/p}$ be a primitive $p$th root of unity and 
	\begin{equation*}
	V=Vandermonde \big( \zeta, \zeta^2, \cdots, \zeta^{p-1} \big)
	\end{equation*}
	Then we have
	\begin{equation*}
	N_q(q-2)=\Big( 1-\frac{1}{q} \Big) \big( (q-1)!-per(V) \big).
	\end{equation*}
	\end{thr}	

	$\it{Proof}$. From $(2)$ it is sufficient to show that
	\begin{equation*}
	\# \big( x_1+\omega x_2+\omega^2 x_3+\cdots +\omega^{q-2} x_{q-1}=0 \big) = \frac{(q-1)!}{q}+\frac{(q-1)per(V)}{q}.
	\end{equation*}
	For the matrix
	\begin{equation*}
	A=\left(
	\begin{array}{cccc}
	X^{\omega} & X^{\omega^2} & \cdots & X^{\omega^{q-1}}\\
	X^{\omega^2} & X^{\omega^3} & \cdots & X^{\omega}\\
	\cdots & \cdots & \cdots & \cdots \\
	X^{\omega^{q-1}} & X^{\omega} & \cdots & X^{\omega^{q-2}}
	\end{array} \right),
	\end{equation*}
	assume that $f(X):=per(A)=\sum_{i=0}^{q-2}c_i X^{\omega^i}+c_{-1}$. If $Tr(\omega^i) \neq 0$, then
	\begin{equation*}
	\sum_{k=0}^{p-1} c_i \xi^{kTr(\omega^i)} = c_i \frac{\xi^{pTr(\omega^i)}-1}{\xi^{Tr(\omega^i)}-1}=0.
	\end{equation*}
	Since $c_0=c_i$ for all $1 \leq i \leq q-2$ from the proof of Theorem 2.1,
	\begin{equation*}
	\sum_{k=0}^{p-1} f\big( \xi^k \big) = p \big( c_{-1}+c_0(p^{r-1}-1) \big). 
	\end{equation*}
	If $k=0$, then $f\big( \xi^k \big)=f(1)=(q-1)!$. If $1 \leq k \leq p-1$, then $f\big( \xi^k \big)=f(\xi)$.\\
	So
	\begin{equation}
	(q-1)!+(p-1)per(V)=p \big( c_{-1}+c_0 (p^{r-1}-1) \big).
	\end{equation}
	From $(3)$ and $(4)$,
	\begin{equation*}
	c_{-1}=\frac{(q-1)!}{q}+\frac{(q-1)per(V)}{q}.
	\end{equation*}
	Then Theorem 2.1 shows that
	\begin{equation*}
	\# \big( x_1+\omega x_2+\omega^2 x_3+\cdots +\omega^{q-2} x_{q-1}=0 \big) = \frac{(q-1)!}{q}+\frac{(q-1)per(V)}{q}.
	\end{equation*}
	Also
	\begin{equation*}
	N_q(q-2)=\Big( 1-\frac{1}{q} \Big) \big( (q-1)!-per(V) \big). ~ \Box
	\end{equation*}

%
%
	\begin{thr}
	For the above matrix, the following holds:
	\begin{equation*}
	per(V)=\sum_{\pi \in \Pi([n])}(-1)^{q-1-t}(d_1-1)!(d_2-1)! \cdots (d_t-1)!C(\pi_1)C(\pi_2) \cdots C(\pi_t)
	\end{equation*}
	where
	\begin{equation*}
	C(\pi_i)=\left \{
	\begin{array}{ll}
		-1, & \alpha \neq 0\\
		q-1, & \alpha = 0
	\end{array} \right.
	\pi_i=\{ i_1, i_2, \cdots, i_k \}, \alpha=\omega^{i_1}+\omega^{i_2}+\cdots+\omega^{i_k}.
	\end{equation*}
	\end{thr}	
	$\it{Proof}$. It is sufficient to calculate $r_{\pi_i}$'s. If
	\begin{equation*}
	\pi_i=\{ i_1, i_2, \cdots, i_k \}, \quad \alpha=\omega^{i_1}+\omega^{i_2}+\cdots+\omega^{i_k},
	\end{equation*}
	then
	\begin{equation*}
	r_{\pi_i}=\xi^{Tr(\alpha)}+\xi^{Tr(\omega \alpha)}+\cdots+\xi^{Tr(\omega^{q-2} \alpha)}.
	\end{equation*}
	For any $\alpha \neq 0, \big\{ \alpha, \omega \alpha, \omega^2 \alpha, \cdots, \omega^{q-2} \alpha \big\} = \textbf{F}_q^{\times}$ and
	\begin{align*}
	r_{\pi_i} & =\big( \xi+\xi^2+\cdots+\xi^{p-1} \big) p^{r-1}+p^{r-1}-1 \\
			& =\frac{1-\xi^{p-1}}{1-\xi} \xi p^{r-1}+p^{r-1}-1=-1.
	\end{align*}
	If $\alpha = 0$, then it is clear that $r_{\pi_i}=q-1$. This completes the proof. ~ $\Box$	
%
%
 %
%
	\section{A bound for the number of permutation polynomials of degree $q-2$}
%
%
	
	\begin{thr}
	The following holds:
	\begin{equation*}
	\Big( 1-\frac{1}{q} \Big) \Big((q-1)!-\frac{1}{q-1}q^{\frac{q}{2}} \Big) \leq N_q(q-2) \leq \Big( 1-\frac{1}{q} \Big) 
	\Big((q-1)!+\frac{1}{q-1}q^{\frac{q}{2}} \Big).
	\end{equation*}
	\end{thr}
	$\it{Proof}$. Let $\lambda_1, \lambda_2, \cdots, \lambda_{q-1}$ be eigenvalues of $V$. Note that $\chi \big( \omega^i \big)=\xi^{Tr(\omega^i)}$ 
	and $\psi_j \big( \omega^k \big)=e^{\frac{2\pi i}{q-1}kj}$ are additive and multiplicative characters in $\textbf{F}_q$ respectively. Hence
	\begin{equation*}
	\lambda_j=\sum_{k=0}^{q-1}\xi^{Tr(\omega^i)} e^{\frac{2\pi i}{q-1}kj} = G(\psi_j, \chi),
	\end{equation*}
	where $G$ is the Gaussian sum.	
	On the other hand, $|\lambda_j|=q^{\frac{1}{2}}$ and this completes the proof. ~ $\Box$
%
%
%
%
	\section{The number of permutation polynomials of degree $d<q-2$}

	Formal polynomials of several variables are defined inductively just as in the usual polynomial theory.
	Let $d$ be a positive integer and $m=q-1-d$. We have the following expression about the number of permutation polynomials of degree $d$.

%
%

	\begin{thr}
	Let $A$ be the $(q-1)\times (q-1)$ matrix $A=\Big( \prod_{l=1}^m X_l^{\omega^{(i-1)l+j}} \Big)_{i, j=\overline{1, q-1}}$ 
	and its permanent be the following;
	\begin{equation*}
	per(A)=\sum_{\substack{(\alpha_1, \cdots, \alpha_m) \in \textnormal{\textbf{F}}_q^m \\ \exists i, \alpha_i\neq 0}}c_{\alpha_1 \alpha_2 \cdots \alpha_m}
	X_1^{\alpha_1}	X_2^{\alpha_2}	\cdots X_m^{\alpha_m}+c.
	\end{equation*}
	Then the number of solutions in $\textnormal{\textbf{F}}_q^{\times}$ of the equations
	\begin{align*}
	& x_1+\omega^{q-d-1} x_2 +\omega^{2(q-d-1)} x_3+\cdots+\omega^{(q-2)(q-d-1)} x_{q-1}=0 \\
	& x_1+\omega^{q-d-2} x_2 +\omega^{2(q-d-2)} x_3+\cdots+\omega^{(q-2)(q-d-2)} x_{q-1}=0 \\
	& \cdots \quad \cdots \quad \cdots \quad \cdots \quad \cdots \\
	& x_1+\omega x_2 +\omega^2 x_3+\cdots+\omega^{q-2} x_{q-1}=0
	\end{align*}
 	such that $x_i \neq 0$ and $x_i \neq x_j$ for $i \neq j$ equals the constant term $c$ of $per(A)$.\\
	\end{thr}
	$\it{Proof}$. From the definition of permanent, every term in $per(A)$ is expressed by
	\begin{equation*}
	\prod_{i=1}^{q-1} \prod_{l=1}^m X_l^{\omega^{(i-1)l+\sigma(i)}}
	\end{equation*}
	for some permutation $\sigma \in S_{q-1}$ and vise versa. Hence
	\begin{equation*}
	per(A)=\sum_{\sigma \in S_{q-1}} \prod_{i=1}^{q-1} \prod_{l=1}^m X_l^{\omega^{(i-1)l+\sigma(i)}}.
	\end{equation*}
	And
	\begin{equation*}
	\prod_{i=1}^{q-1} \prod_{l=1}^m X_l^{\omega^{(i-1)l+\sigma(i)}} = \prod_{l=1}^m \prod_{i=1}^{q-1} X_l^{\omega^{(i-1)l} \omega^{\sigma(i)}} 
	= \prod_{l=1}^m  X_l^{\sum_{i=1}^{q-1}\omega^{(i-1)l}\omega^{\sigma(i)}}.
	\end{equation*}\\
	On the other hand, it is clear that $X_l^{\sum_{i=1}^{q-1}\omega^{(i-1)l}\omega^{\sigma(i)}}=1$ if and only if
	\begin{equation*}
	\omega^{\sigma(1)}+\omega^l \omega^{\sigma(2)}+\omega^{2l} \omega^{\sigma(3)}+ \cdots +\omega^{(q-2)l} \omega^{\sigma(q-1)}=0.
	\end{equation*}
	So
 	\begin{align*}
	& x_1+\omega^m x_2 +\omega^{2m} x_3+\cdots+\omega^{(q-2)m} x_{q-1}=0 \\
	& x_1+\omega^{m-1} x_2 +\omega^{2(m-1)} x_3+\cdots+\omega^{(q-2)(m-1)} x_{q-1}=0 \\
	& \cdots \quad \cdots \quad \cdots \quad \cdots \quad \cdots \\
	& x_1+\omega x_2 +\omega^2 x_3+\cdots+\omega^{q-2} x_{q-1}=0
	\end{align*}
	for $x_i \in \textbf{F}_q^{\times}$, which are different from each other, if and only if there is a permutation $\sigma \in S_{q-1}$ 
	such that $\prod_{l=1}^m  X_l^{\sum_{i=1}^{q-1}\omega^{(i-1)l}\omega^{\sigma(i)}}=1$ and $x_i = \omega^{\sigma (i)}$ 
	for all integers $i \in \{1, 2, \cdots, q-1 \}$. Thus $G_q (d)=c$ and this completes the proof. ~ $\Box$
%
%
%
%
	\section{Conclusion and further study}
	In this paper we have obtained a formula for the number of permutation polynomials of degree $\leq q-2$ and a bound for the number of 
	permutation polynomials of degree $q-2$ over $\textbf{F}_q$ by generalizing the results in \cite{das}. 
	The problem of finding a new better bound for the number of permutation polynomials of degree $d-2$ will be worth exploring.\\

	{\bf Acknowledgement}. We would like to thank anonymous referees for their valuable comments and suggestions.

 	\end{document}